\documentclass[fleqn,11pt]{article}
\usepackage{graphicx}
\usepackage{amscd,amssymb}
\usepackage{amsmath}

\title{A  finite time blowup result for quadratic ODE's}
\author{Dennis Sullivan*}
\date{}

\begin{document}

\maketitle

 ....dedicated to Mauricio Peixoto who abstracted the practical theory of ODE's and to David Rand who applied the subsequent powerful abstract theory to practical problems.

The famous Euler ODE of incompressible frictionless fluid dynamics expressed in terms of the variable $X$ = vorticity has the following algebraic form:
 The underlying space can be viewed as the (infinite dimensional) vector space $V$ of exact two forms on a closed Riemannian manifold.
The evolution of  the exact two form = vorticity is described by an ODE
 $dX/dt = Q(X)$ where $Q$ is a  homogeneus quadratic function on V, namely one whose deviation from being linear $Q(X+Y)$ - $Q(X)$ $-Q(Y)$ is a symmetric bilinear form on $V$.
We make a few comments about finite time blowup with given initial conditions for such algebraic, specifically homogeneous quadratic, ODE's.

1) Of course the most simple example on the real line $dx/dt$ =$- x.x$ with initial condition a at $t=0$ has solution $1/(t-1/a)$. This blows up at the critical time $t =1/a$.

2) The same calculation works for $dX/dt$ =$- X.X$ where $X$ is a linear operator , so $V= End(W)$ for some other linear space $W$. If $A$ is the operator at time zero
 Then $X(t)$ = $Id/( t(Id)- Id/A)$ is the solution. So $X(t)$ blows up at a finite time iff the spectrum of $A$ contains a real number. So for an open set of initial conditions there is a solution for all time and for another open set there is finite time blowup.

3) One may hope  to find  some structure like this in the Euler fluid equation referred to above that would prevent finite time blowup for a large set of initial conditions.

Thus we ask the following questions :
1) for a finite dimensional vector space V how likely is it in the variable $Q$ for the quadratic ODE  $dx/dt$ = $Q(X)$ to have finite time blow up for some initial condition $A$.
2) fixing $Q$ how likely is it in the variable $A$ to have finite time blowup.

The following theorem answers question 1)

Theorem: If $V$ is any finite dimensional vector space, then outside a proper algebraic subvariety of quadratic functions $Q$, the ODE   $dX/dt$ = $Q(X)$
 exhibits finite time blowup for some initial condition.

Proof: (Ofer Gabber) The condition that there exists a non zero vector $Y$ so that $Q(Y)$ = $0$ defines a proper algebraic subvariety in the space of quadratic mappings from $V$ to $V$.
 Outside this subvariety $Q$ defines by rescaling a map from $S$, the sphere of directions in $V$,  to itself. This mapping agrees on antipodal points because $Q$ is quadratic.
Thus this mapping  from $S$ to $S$ has even topological degree.Then the lefschetz number of this mapping is non zero and this mapping has a fixed point. This means the original $Q$ keeps a line invariant. By a linear change of variable the ODE restricted to this line becomes example 1) which  has finite time blowup. QED.

We have not studied question 2)  in general which in example 2) is  interesting.

* The theorem came out of a two day discussion at IHES in the early 90's with Ofer Gabber where the author made up the question 1) and Ofer came up with the proof of the Theorem.

\end{document}